\documentclass[11pt]{article}
\author{M.~E.~Kazarian\thanks{Steklov Mathematical Institute
RAS and the Poncelet Laboratory, Independent University of Moscow.
Partly supported by NWO-RFBR grant 047.011.2004.026 (RFBR05-02-89000 NWOa),
RFBR grant 04-01-00762, ANR GIMP ANR-05-BLAN-0029-01.}, S.~K.~Lando
\thanks{Institute for System Research RAS and the Poncelet Laboratory, Independent
University of Moscow. Partly supported by the grants
 ACI-NIM-2004-243
(Noeuds et tresses), NWO-RFBR 047.011.2004.026
(RFBR 05-02-89000-NWOa), 
RFBR 02-01-22004-CNRSa,
05-01-01012-a and ANR GIMP ANR-05-BLAN-0029-01.}}

\title{An algebro-geometric proof of Witten's conjecture}
\date{January 31, 2006}

\usepackage{amssymb,amsmath}
\usepackage{theorem}

\def\C{{\mathbb C}}

\def\Z{{\mathbb Z}}

\def\cF{{\cal F}}

\def\cL{{\cal L}}

\def\ocM{{\overline{{\cal M}}}}

\def\st{{\rm st}}

\def\Aut{{\rm Aut}}

\def\p{{\partial}}

\def\dim{{\rm dim}}

\newtheorem{theorem}{Theorem}[section]

\newtheorem{corollary}[theorem]{Corollary}

{\theorembodyfont{\rmfamily}    

}

\begin{document}

\maketitle

\begin{abstract}
\noindent
We present a new proof of Witten's conjecture. The proof is based
on the analysis of the relationship between intersection indices
on moduli spaces of complex curves and Hurwitz numbers enumerating
ramified coverings of the $2$-sphere.
\end{abstract}

\section{Introduction}
Let~$\ocM_{g;n}$ denote the Knudsen--Deligne--Mumford 
moduli space of genus~$g$ stable complex curves
with~$n$ marked points~\cite{DM}. For each $i\in~\{1,\dots,n\}$,
consider the line bundle ~$\cL_i$ over~$\ocM_{g;n}$ 
whose fiber over a point $(C;x_1,\dots,x_n)\in\ocM_{g;n}$
is the cotangent line to the curve~$C$
at the marked point~$x_i$. Let~$\psi_i\in H^2(\ocM_{g;n})$
denote the first Chern class of this line bundle,
$\psi_i=c_1(\cL_i)$. Consider the generating function
\begin{equation}\label{Wit}
F(t_0,t_1,\dots)=\sum\langle\tau_{d_1}\dots\tau_{d_n}\rangle
\frac{t_{d_1}\dots t_{d_n}}{|\Aut(d_1,\dots,d_n)|}
\end{equation}
for the intersection numbers of these classes,
$$
\langle\tau_{d_1}\dots\tau_{d_n}\rangle=
\int_{\ocM{g;n}}\psi_1^{d_1}\dots \psi_n^{d_n},
$$
where the genus~$g$ is uniquely determined from the identity
$$
d_1+\dots+d_n=\dim~\ocM_{g;n}=3g-3+n.
$$
The first few terms of this function are
\begin{eqnarray*}
F&=& \frac1{24}t_1+\frac1{6}t_0^3+\frac1{48}t_1^2
+\frac1{24}t_0t_2+\frac1{6}t_0^3t_1+
\frac1{1152}t_4+\frac1{72}t_1^3+\frac1{12}t_0t_1t_2
+\frac1{48}t_0^2t_3\\
&&+\frac1{6}t_0^3t_1^2+
\frac1{24}t_0^4t_2+\frac{29}{5760}t_2t_3+\frac1{384}t_1t_4+
\frac1{1152}t_0t_5+\dots
\end{eqnarray*}

The celebrated Witten conjecture asserts that

{\it the
second derivative $U=\p^2F/\p t_0^2$ of the generating function~$F$
satisfies the KdV equation}
\begin{equation}\label{KdV}
\frac{\p U}{\p t_1}=U\frac{\p U}{\p t_0}+\frac1{12}\frac{\p^3 U}{\p t_0^3}.
\end{equation}

The motivation for this conjecture can be found in~\cite{W1},
and for a detailed exposition suitable for a mathematically-minded
reader see, e.g.,~\cite{LZ}. Note that the KdV equation can be interpreted
as a reccurence formula allowing one to calculate all the intersection
indices provided ``initial conditions'' are given.
Witten has shown that the function~$F$ satisifes the so-called string
and dilaton equations, which together with the KdV equation generate
the whole KdV hierarchy and provide necessary initial conditions.

A number of proofs of the Witten conjecture are known,
but all of them exploit techniques that do not seem to be
intrinsically related to the initial problem: Kontsevich's
proof~\cite{K} makes use of Jenkins--Strebel differentials and matrix
integrals, the proof due to Okounkov and Pandharipande~\cite{OP},
which starts with Hurwitz numbers, also involves matrix
integrals and graphs on surfaces study, as well as asymptotic analysis,
and, finally, Mirzakhani's proof~\cite{M}
is based on the Riemannian geometry properties of moduli spaces.
The goal of this paper is to present a new proof using purely
algebro-geometric techniques. Similarly to the proof due to 
Okounkov and Pandharipande, we start with Hurwitz numbers,
but then we follow a different line.

Hurwitz numbers enumerate ramified coverings of the $2$-sphere
with prescribed ramification points and ramification types
over these points. We deal only with ramified coverings
whose ramification type is simple over each ramification
point but one. Our proof is based on the following, now
well-known, properties of these numbers:
\begin{itemize}
\item the ELSV formula~\cite{ELSV1,ELSV2}, see also~\cite{GV}, \cite{L},
relating Hurwitz numbers to the intersection theory
on moduli spaces;
\item the relationship between Hurwitz 
numbers and integrable hierarchies
conjectured by Pandharipande~\cite{P} and proved,
in a stronger form, by Okounkov~\cite{O}.
\end{itemize}
Using the ELSV formula we express the intersection indices of
the $\psi$-classes in terms of Hurwitz numbers. The
partial differential equations governing the generating series for
Hurwitz numbers then lead to the KdV equation for the intersection indices.
Note that the existence of such a proof has been predicted in~\cite{GJV}.
One of the main features of the proof consists in the fact
that known effective algorithms for computing
the Hurwitz numbers, which are relatively
simple combinatorial objects, lead to an independent tool for computing
the intersection indices. 
We describe the properties of the Hurwitz numbers
in detail and deduce Witten's conjecture
from them in Sec.~2. Section~3 is devoted to a discussion of the proof.

We express our gratitude to the participants 
of our seminar at the Independent
University of Moscow
for useful discussions. Our personal gratitude
is due to S.~Shadrin and D.~Zvonkine, whose papers contain,
in an implicit form, the idea of inverting the ELSV formula.
Special thanks are due to D.~Zvonkine for careful 
proofreading of the first versions of the paper.
The final version of this paper
has been written at the Max-Planck Institut f\"ur Mathematik, 
Bonn (MPI preprint 05-55),
to whom we are grateful for hospitality.

\section{Proof}

\subsection{Hurwitz numbers}
Fix a sequence $b_1,\dots,b_n$ of positive integers. 
Consider ramified coverings of the sphere~$S^2$
by compact oriented two-dimensional surfaces
of genus~$g$ with ramification type $(b_1,\dots,b_n)$
over one point, and the simplest possible ramification
type $(2,1,1,\dots,1)$ over all other points of ramification.
According to the Riemann--Hurwitz formula, the total
number~$m$ of these points of simple ramification is
\begin{equation}\label{dim}
m=2g-2+n+B,
\end{equation}
where $B=b_1+\dots+b_n$ is the degree of the covering.
If we fix the ramification points in the target sphere,
then the number of topologically distinct ramified
coverings of this type becomes finite, and we denote
by $h_{g;b_1,\dots,b_n}$ the number of these coverings,
with marked preimages of the point of degenerate ramification, 
counted with the weight inverse to the order of the
automorphism group of the covering. These numbers
are called {\it Hurwitz numbers}.

The ELSV formula~\cite{ELSV1,ELSV2} expresses the
Hurwitz numbers in terms of Hodge integrals over the moduli spaces
of stable complex curves:
\begin{equation}\label{ELSV}
h_{g;b_1,\dots,b_n}=m!
\prod_{i=1}^n\frac{b_i^{b_i}}{b_i!}
\int_{\ocM_{g;n}}\frac{1-\lambda_1+\lambda_2-\dots\pm\lambda_g}
{(1-b_1\psi_1)\dots(1-b_n\psi_n)}
\end{equation}
for $g>0,n\ge1$ or $g=0,n\ge3$.
The numerator of the integrand is the 
total Chern class of the vector bundle over $\ocM_{g;n}$
dual to the Hodge bundle (whose fiber is 
the $g$-dimensional vector space
of holomorphic differentials on the curve), 
$\lambda_i\in H^{2i}(\ocM_{g;n})$.
The integral in~(\ref{ELSV}) is understood as the result
of expanding the fraction as a power series, with further
selection of monomials whose degree coincides with the
dimension of the base (there are finitely many of them)
and integration of each of these monomials.

The integral~(\ref{ELSV}) is a sum of intersection indices
of both $\psi$- and $\lambda$-classes. In~\cite{OP}, the
$\lambda$-classes are eliminated by considering asymptotics
of integrals of this kind. In contrast, in the present paper, the
exclusion of the $\lambda$-classes is based on simple
combinatorial considerations originating in~\cite{Sh,Z},
see Sec.~2.2 below.

Now consider the following exponential generating function
for the Hurwitz numbers:
$$
H(\beta;p_1,p_2,\dots)=
\sum h_{g;b_1,\dots,b_n} \frac{p_{b_1}\dots p_{b_n}}{n!} \frac{\beta^m}{m!},
$$
where the summation is taken over all finite sequences
$b_1,\dots,b_n$ of positive integers and all nonnegative
values of~$g$, with~$m$ given by Eq.~(\ref{dim}). 
According to Okounkov~\cite{O}, the exponent $e^H$ of
this generating function is a $\tau$-function for the KP-hierarchy.
In fact, Okounkov proved a much stronger theorem stating
that the generating function for double Hurwitz numbers
(those having degenerate ramification over two points rather than one)
satisfies the Toda lattice equations. We do not need this statement
in such generality, and in Sec.~3 we discuss a simple proof
of the fact we really need.
(Various kinds of Hurwitz numbers have been since long known to lead
to solutions of integrable hierarchies, 
but we were unable to trace exact statements and origins.)
The function $e^H$ being a $\tau$-function for the KP hierarchy
means, in particular, that $H$ satisfies the {\it KP-equation}
\begin{equation}\label{KP}
\frac{\p^2 H}{\p p^2_2}=
\frac{\p^2H}{\p p_3\p p_1}-\frac12\left(\frac{\p^2 H}{\p p_1^2}\right)^2
-\frac1{12}\frac{\p^4 H}{\p p^4_1}.
\end{equation}

We use this equation below to deduce
from it the KdV equation for the function~$U$.

\subsection{Expressing intersection indices of $\psi$-classes
via Hurwitz numbers}

Obviously, for each nonegative integer~$d$ there exist constants~$c_b^d$,
$b=1,\dots,d+1$, such that
\begin{equation}\label{lk}
 \sum_{b=1}^{d+1}\frac{c_b^d}{1-b\psi}
 =\psi^d+O(\psi^{d+1}),
\end{equation}
and these constants are uniquely determined by this requirement.
They are given by the formula
$$
c_b^d=\frac{(-1)^{d-b+1}}{(d-b+1)!(b-1)!}.
$$
Indeed, we need to prove that the first $d-1$ derivatives in~$\psi$ of
the linear combination
$$
\sum_{b=1}^{d+1}\frac{c_b^d}{1-b\psi}
$$
vanish at~$0$, while the $d$~th derivative is~$d!$. The~$i$~th derivative of
this linear combination evaluated at~$\psi=0$ is 
$$
(-1)^{i+1}\frac1{i!}\left(\binom{d}{0}1^i-\binom{d}{1}2^i+\binom{d}{2}3^i-\dots\pm
\binom{d}{d}d^i\right).
$$
The expression in brackets coincides with the result of applying the 
$i$~th iteration of the operator $xd/dx$ to the polynomial $(1-x)^d$
and evaluating at~$x=1$, which is~$0$ for $0\le i<d$ and $(-1)^d d!$ for $i=d$.

Multiplying identities~(\ref{lk}) for different $d$ we obtain
the following equality:
$$
\sum_{b_1=1}^{d_1+1}\cdots\sum_{b_n=1}^{d_n+1}
 \frac{c_{b_1}^{d_1}\dots c_{b_n}^{d_n}}
 {(1-b_1\,\psi_1)\dots(1-b_n\,\psi_n)}
 =\prod_{i=1}^n\psi_i^{d_i}+\dots,
$$
where dots on the right-hand side denote cohomology classes
of degree greater than $d_1+\dots+d_n$. This means, in particular,
that for $d_1+\dots+d_n=3g-3+n$ the linear combination
$$
\sum_{b_1=1}^{d_1+1}\cdots\sum_{b_n=1}^{d_n+1}
c_{b_1}^{d_1}\dots c_{b_n}^{d_n}
\int_{\ocM_{g;n}}
\frac{1-\lambda_1+\dots\pm\lambda_g}{(1-b_1\psi_1)\dots(1-b_n\psi_n)}
$$
is simply $\langle\tau_{d_1}\dots\tau_{d_n}\rangle$, because the integral
of the terms of higher degree vanishes.
Taking into account the coefficient of the integral 
in Eq.~\eqref{ELSV}, we obtain the following explicit identity.

\begin{theorem}
For any sequence of non-negative integers $d_1,\dots,d_n$ we have
$$
\langle\tau_{d_1}\dots\tau_{d_n}\rangle=
\sum_{b_1=1}^{d_1+1}\dots
\sum_{b_n=1}^{d_n+1}\left(\frac1{m!}
\prod_{i=1}^n\frac{(-1)^{d_i+1-b_i}}{(d_i+1-b_i)!b_i^{b_i-1}}\right)
h_{g;b_1,\dots,b_n},
$$
where $g$ is determined by the left-hand side, $\sum d_i=3g-3+n$,
and $m=2g-2+B+n$.
\end{theorem}

It is convenient to reformulate the statement of the theorem in terms
of generating functions. Decompose the generating function~$H$ into the sum
\begin{equation}\label{gen}
H=H_{0;1}+H_{0;2}+H_\st,
\end{equation}
where the stable part $H_\st$ contains all the monomials whose coefficients
are given by the ELSV formula~(\ref{ELSV}), and $H_{0;1}$
and $H_{0;2}$ are the generating functions
for the numbers of ramified coverings of the sphere by the sphere
with~$1$ (``polynomial'') and~$2$ (``trigonometric polynomial'')
preimages over the distinguished ramification point, respectively.
The latter generating functions are known since Hurwitz~\cite{H1,H2}, 
see also~\cite{A}:
\begin{eqnarray*}
H_{0;1}&=&\sum_{b=1}^\infty h_{0;b} p_b \frac{\beta^{b-1}}{(b-1)!} = 
\sum_{b=1}^\infty \frac{b^{b-2}}{b!} p_b \beta^{b-1} \\
H_{0;2}&=&\frac12\sum_{b_1,b_2=1}^\infty 
h_{0;b_1,b_2} p_{b_1}p_{b_2}\frac{\beta^{b_1+b_2}}{(b_1+b_2)!} =\frac12
\sum_{b_1,b_2=1}^\infty 
\frac{b_1^{b_1}b_2^{b_2}}{(b_1+b_2)b_1!b_2!}
p_{b_1}p_{b_2}\beta^{b_1+b_2} 
\end{eqnarray*}
(note that this case can also be considered as been covered
by the ELSV formula, but with the moduli {\it spaces} $\ocM_{g;n}$
replaced by the moduli {\it stacks} $\ocM_{0;1}$, $\ocM_{0;2}$).
In fact, we are going to use below not the precise formulas for
$H_{0;1},H_{0;2}$, but the fact that they contain only terms
of degree at most~$2$ in~$p_i$, which yields
\begin{eqnarray*}
\frac{\p^2}{\p p_1^2}(H_{0;1}+H_{0;2})&=&\frac12\beta^2,\\
\frac{\p^2}{\p p_2^2}(H_{0;1}+H_{0;2})&=&\beta^4,\\
\frac{\p^2}{\p p_1\p p_3}(H_{0;1}+H_{0;2})&=&\frac98\beta^4.
\end{eqnarray*}
Denote by $G_\st=G_\st(\beta;t_0,t_1,\dots)$ 
the result of the following change of variables in the series $H_\st$:
\begin{equation}\label{subs}
p_b=\sum_{d=b-1}^\infty\frac{(-1)^{d-b+1}}{(d-b+1)!b^{b-1}}
\beta^{-b-\frac{2d+1}3}t_d.
\end{equation}
The result of this substitution is a series in $t_0,t_1,\dots$ whose
coefficients are formal Laurent expansions in $\beta^{2/3}$. Indeed,
the powers of~$\beta$ in the contribution of a monomial 
$p_{b_1}\dots p_{b_n}$ to the expansion have the form
$$
(b_1+\dots+b_n+2g-2+n)+\sum_{i=1}^n\left(-b_i-\frac{2d_i+1}3\right)=
\frac23\left(3g-3+n-\sum_{i=1}^nd_i\right),
$$
hence become even integers when multiplied by~$3$. On the other hand,
the powers of~$\beta$ at each~$t_d$ are bounded from below,
because~$t_d$ enters only the expansions for $p_1,\dots,p_{d+1}$.

\begin{theorem}\label{II}
{\rm1.} The series $G_\st$ contains no terms with negative powers of $\beta$.

{\rm2.} The free term in $\beta$, $G_\st|_\beta=0$ 
{\rm(}which is correctly defined due to the first statement{\rm)}
coincides with 
the generating function~$F$ for the intersection numbers given by 
Eq.~{\rm(\ref{Wit})},
$$
F(t_0,t_1,\dots)=G_\st(0;t_0,t_1,\dots).
$$
\end{theorem}

{\bf Proof.} Collect together the terms of the 
series $H_\st$ corresponding to given values $n,b_1,\dots,b_n$,
for all~$g$, and set
$$
H_{b_1,\dots,b_n}=\sum_{m=B+2g-2+n} 
h_{g;b_1,\dots,b_n}\frac{\beta^m}{m!}.
$$
Then the ELSV formula~(\ref{ELSV})
can be conveniently rewritten as
$$ 
 H_{b_1,\dots,b_n}=\beta^{B+n/3}
\prod_{i=1}^n\frac{b_i^{b_i}}{b_i!}
 \left\langle\frac{1-\beta^{2/3}\lambda_1+\beta^{4/3}\lambda_2-
\beta^{6/3}\lambda_3+\dots}
{(1-b_1\beta^{2/3}\psi_1)\cdots(1-b_n\beta^{2/3}\psi_n)}\right\rangle,$$
where we understand the numerator as the formal sum and
the angle brackets mean integration of each monomial
over the space~$\ocM_{g;n}$ whose dimension coincides
with the degree of the monomial. Indeed, 
consider a summand containing the integral of $\psi_1^{d_1}\dots\psi_n^{d_n}\lambda_j$
 in the expansion on the right-hand 
side. The genus $g$ corresponding to the domain of integration $\ocM_{g;n}$ is computed from the equality
$$\sum d_i+j=3g-3+n.$$
This summand contributes to a term of degree $m$ in $\beta$ iff 
the relation 
$$
m=B+\frac n3+\frac23\sum_{i=1}^n d_i+\frac23j=2g-2+n+B,
$$
which is exactly the relation between $g$ and $m$ in the definition of the series $H_{b_1,\dots,b_n}$, holds.

Now, the explicit form of the change
of variables~(\ref{subs}) implies
that the coefficient $G_{d_1,\dots,d_n}$ of the monomial $t_{d_1}\dots t_{d_n}/|\Aut(d_1,\dots,d_n)|$ 
in~$G_\st$ is equal to
\begin{align*}
G_{d_1,\dots,d_n}&= \sum_{b_1=1}^{d_1+1}\dots\sum_{b_n=1}^{d_n+1} \left(
 \prod_{i=1}^n
\tfrac{(-1)^{d_i-b_i+1}\beta^{-b_i-\frac{2d_i+1}3}}
{(d_i-b_i+1)!b_i^{b_i-1}}\right)H_{b_1,\dots,b_n}\\
&=\beta^{-\frac23\sum d_i} \sum_{b_1=1}^{d_1+1}\dots\sum_{b_n=1}^{d_n+1} \left\langle
\frac{c_{b_1}^{d_1}\dots c_{b_n}^{d_n}(1-\beta^{2/3}\lambda_1+\beta^{4/3}\lambda_2-
\beta^{6/3}\lambda_3+\dots)}
{(1-b_1\beta^{2/3}\psi_1)\cdots(1-b_n\beta^{2/3}\psi_n)}\right\rangle\\
&=\beta^{-\frac23\sum d_i}\langle(\prod_{i=1}^n(\beta^{2/3}\psi_i)^{d_i} +
\dots) (1-\beta^{2/3}\lambda_1+\beta^{4/3}\lambda_2-
\dots)\rangle\\
&=\langle\prod_{i=1}^n\psi_i^{d_i}\rangle+\dots,
\end{align*}
where dots in the last two lines denote terms of higher degree in~$\beta$.

Theorem~\ref{II} is proved.

\subsection{Reduction to the KdV equation}
By Theorem~\ref{II}, $G_\st$
is a power series in $\beta^{2/3},t_0,t_1,\dots$
whose coefficient $G|_{\beta=0}$ of~$\beta^0$ coincides with the
function $F$. By definition, $G_\st$ is the result of the substitution~(\ref{subs}) 
to the series
$H_\st=H-H_{0;1}-H_{0;2}$.
Equation~(\ref{KP}) together with the explicit formulas for $H_{0;1}$ and $H_{0;2}$
imply that $H_\st$ satisfies the equation
$$
\frac{\p^2 H_\st}{\p p_2^2}=
\frac{\p^2 H_\st}{\p p_3\p p_1}-
\frac12\left(\frac{\p^2 H_\st}{\p p_1^2}\right)^2
-\frac12\beta^2\frac{\p^2 H_\st}{\p p_1^2}
-\frac1{12}\frac{\p^4 H_\st}{\p p_1^4}.
$$

The change of variables~(\ref{subs}) results in the following
change of partial derivatives:
\begin{eqnarray*}
\frac{\p }{\p p_1}&=&\beta^{4/3}\frac{\p }{\p t_0};\\
\frac{\p }{\p p_2}&=&2\beta^{9/3}\frac{\p }{\p t_1}+2\beta^{7/3}\frac{\p}{\p t_0}\\
\frac{\p }{\p p_3}&=&9\beta^{14/3}\frac{\p }{\p t_2}
+9\beta^{12/3}\frac{\p }{\p t_1}+\frac92\beta^{10/3}\frac{\p}{\p t_0}.
\end{eqnarray*}
After substituting this into the equation above
and dividing the result by~$\beta^{16/3}$, we rewrite it as
\begin{equation}\label{preKdV}
\left(\frac{\p^2 G_\st}{\p t_1\p t_0}-\frac12\left(\frac{\p^2 G_\st}{\p t_0^2}\right)^2
-\frac1{12}\frac{\p^4 G_\st}{\p t_0^4}\right)
+\beta^{2/3}\left(9\frac{\p^2 G_\st}{\p t_0\p t_2}-4\frac{\p^2 G_\st}{\p t^4_0}
\right)=0.
\end{equation}
The coefficient of~$\beta^0$ in~(\ref{preKdV}) 
has the form
$$
\frac{\p^2 F}{\p t_1\p t_0}-\frac12\left(\frac{\p^2 F}{\p t_0^2}\right)^2-
\frac1{12}\frac{\p^4 F}{\p t_0^4}=0.
$$
Differentiating the last equation over~$t_0$ we obtain the
desired KdV equation~(\ref{KdV}) for the function $U=\p^2 F/\p^2 t_0^2$.
This completes the proof of Witten's conjecture.

\section{Odds and ends}
Trying to make the present paper
more self-contained, 
we discuss in this section the Hirota bilinear equations
and the KP hierarchy, as well as an explicit presentation
of the function $e^H$ as a $\tau$-function
for the KP hierarchy.
Although specialists in integrable hierarchies are well
aware of these facts, it is not an easy task to find
their compact and readable exposition. We refer the reader
to~\cite{S}, \cite{MJT}, \cite{DKJM}, and~\cite{SW} for a description
of the relationship between integrable hierarchies and
the geometry of semi-infinite Grassmannian.
The KP hierarchy is a system of partial differential equations
for an unknown function~$H$. The KP equation~(\ref{KP})
is the first equation in this system. Similarly to the case of
the KdV equation, the expansion of a solution to the KP equation
can be reconstructed from ``initial conditions''.
The exponent $\tau=e^H$ of a solution~$H$ to the KP hierarchy is called
a $\tau$-{\it function} of the hierarchy. The equations of the
KP hierarchy rewritten for $\tau$-functions also are partial
differential equations; they are called
the Hirota equations. They possess a nice property of being
quadratic with respect to~$\tau$.

\subsection{Semi-infinite Grassmannian, Hirota--Pl\"ucker\\
bilinear equations, and integrable hierarchies}

Define the {\it charge zero Fock space}~$\cF$
as the completion of the infinite dimensional coordinate
vector space over~$\C$ whose basic elements~$s_\lambda$ are
labeled by partitions,
$$
\cF=\overline{\bigoplus\C s_\lambda}
$$
Recall that a {\it partition} is a nonincreasing
sequence of integers $\lambda=(\lambda_1,\lambda_2,\dots)$,
$\lambda_1\ge\lambda_2\ge\dots\ge0$, having finitely many nonzero terms.
Elements of~$\cF$ are infinite formal linear combinations of the 
vectors~$s_\lambda$. We shall use the following two
realizations of the Fock space.

(1) The space~$\cF$ can be identified with the space 
$\cF=\C[[p_1,p_2,\dots]]$ of formal power series
in infinitely many variables $p_1,p_2,\dots$
by setting~$s_\lambda$ to be the corresponding
Schur function. The {\it Schur function}
corresponding to a one-part partition is defined
by the expansion
$$
s_0+s_1z+s_2z^2+s_3z^3+s_4z^4+\dots
=e^{p_1z+p_2\frac{z^2}2+ p_3\frac{z^3}3+\dots},
$$
and for a general partition~$\lambda$ it is given by the determinant
\begin{equation}\label{Schur}
s_\lambda=\det~||s_{\lambda_j-j+i}||.
\end{equation}
The indices~$i,j$ here run over the set $\{1,2,\dots,n\}$
for $n$ large enough, and since $\lambda_i=0$ for~$i$
sufficiently large, the determinant, whence~$s_\lambda$,
is independent of~$n$. Here are a few first Schur polynomials:
$$
s_0=1,\qquad s_1=p_1,\qquad s_2=\frac12(p_1^2+p_2),\qquad
s_3=\frac16(p_1^3+3p_1p_2+2p_3),
$$
$$
s_{1,1}=\frac12(p_1^2-p_2), \qquad
s_{2,1}=\frac13(p_1^3-p_3), \qquad
s_{1,1,1}=\frac16(p_1^3-3p_1p_2+2p_3).
$$

(2) Let $V=\C[z,z^{-1}]$ be the ring of Laurent polynomials in~$z$,
which we treat as the vector space with the basis $z^i$, $i\in\Z$.
Identify~$\cF$ with the {\it semi-infinite wedge space}
$\cF=\Lambda^{\infty/2}V$ freely spanned by the formal infinite wedge
products of the form
$$
s_\lambda=z^{k_1}\wedge z^{k_2}\wedge z^{k_3}\wedge\dots,\qquad
k_i=i-\lambda_i,
$$
for all partions~$\lambda$. Sequences~$k_i$ appearing on the right-hand
side can be characterized as arbitrary strictly increasing sequences
of integers satisfying $k_i=i$ for~$i$ large enough.

The theory of the KP hierarchy can be summarized as follows.

The projectivization $P\cF=P\Lambda^{\infty/2}V$ is the ambient
space of the standard Pl\"ucker embedding $Gr\hookrightarrow P\cF$, where
$Gr=Gr_{\infty/2}(V)$ is the Grassmannian of ``half-infinite
dimensional subspaces'', often referred to as the {\it Sato Grassmannian}.
By definition, the elements of~$Gr$
are subspaces spanned by linearly independent vectors
$\varphi_1,\varphi_2,\varphi_3,\dots$ in (the formal completion of)~$V$
such that for~$i$ large enough we have $\varphi_i=z^i+\dots$,
where dots denote terms of lower order in~$z$.
Such a vector space can be interpreted as the wedge product
\begin{equation}\label{wedge}
\tau=\varphi_1\wedge\varphi_2\wedge\varphi_3\wedge\dots
\end{equation}
Indeed, another choice of a basis does not affect this wedge product,
up to a scalar factor.
If all the functions~$\varphi_i$ are Laurent polynomials
and $\varphi_i$ is simply the monomial $z^i$ for all~$i$ sufficiently large,
then the wedge product~$\tau$ can be represented, after
expanding the brackets, as a {\it finite} linear combination
$\tau=\sum_\lambda c_\lambda s_\lambda$, $c_\lambda\in\C$. 
If the functions~$\varphi_i$ contain infinitely many terms,
then the function $\tau$ is a formal linear combination
of~$s_\lambda$ and can be obtained in the following way:

{\it when expanding the brackets in the infinite wedge 
product~{\rm(\ref{wedge})} pick one monomial summand in each~$\varphi_i$
in such a way that this summand is~$z^i$ for all but
finitely many indices~$i$ and do this in all possible ways.}

More explicitly, if $\varphi_i=\sum_{j\in\Z}a_{i,j}z^j$,
then
\begin{equation}\label{tau}
\tau=\sum_\lambda\det||a_{i,j-\lambda_i}||_{i,j\ge1}s_\lambda
=\det||\sum_{k\in\Z}a_{i,k}s_{j-k}||_{i,j\ge1}.
\end{equation}

\begin{theorem}[\cite{S,DKJM,SW}]
A {\rm(}non-zero{\rm)} function~$\tau$ 
is a $\tau$-function for the KP hierarchy
iff the corresponding point $[\tau]\in P\cF$ belongs to the
Grassmannian $Gr\subset P\cF$.
\end{theorem}

In particular, each Schur polynomial~$s_\lambda$
and any linear combination of the Schur functions~$s_i$
corresponding to one-part partitions
produces a solution to the KP equation~(\ref{KP}).

The image of the Pl\"ucker embedding of the Grassmannian is known
to be given by a system of quadratic equations. In our case,
these algebraic equations for the Taylor coefficients of the
function~$\tau$ take the form of partial differential equations
for this function, known
as the {\it bilinear Hirota equations}.

\subsection{Formulas for the generating function for Hurwitz numbers}
The exponent~$e^H$ of the generating function for the Hurwitz numbers
is nothing but the generating function for the numbers of
ramified coverings of the $2$-sphere by {\it all}, not
necessarily connected, compact oriented surfaces of Euler characteristic~$2-2g$.
Take such a covering and let a point of simple ramification
in the target sphere tend to the point of degenerate ramification.
Then  one can express the number of such coverings as a linear
combination of the numbers of similar coverings
with fewer points of simple ramification. This reccurence
relation (the ``cut-and-join equation'' of~\cite{GJ}) 
expressed in terms of generating functions
reads as follows:
\begin{equation}\label{cj}
\frac{\p e^H}{\p\beta}=\frac12\sum_{i,j=1}^\infty
\left((i+j)p_ip_j\frac\p{\p p_{i+j}}+ijp_{i+j}\frac{\p^2}{\p p_i\p p_j}\right)e^H=Ae^H,
\end{equation}
where we denote by~$A$ the differential operator on the right-hand side.
For an algebro-geometric interpretation of the cut-and-join equation,
see~\cite{Sh}. 

Equation~(\ref{cj}) can be solved explicitly. The operator~$A$
acts linearly on the space of weighted homogeneous polynomials
in the variables~$p_1,p_2,\dots$, with the weight of the variable~$p_i$
equal to~$i$. Moreover, it preserves the weighted degree
of the polynomials, whence can be split into a direct sum of finite
dimensional linear operators. The Schur functions~$s_\lambda$ are eigenvectors of~$A$,
and they form a complete set of eigenvectors.

Denote by~$f(\lambda)$ the eigenvalue of the eignevector~$s_\lambda$.
It can be checked that
$$
f(\lambda)=\frac12\sum_{i=1}^\infty\lambda_i(\lambda_i-2i+1).
$$
It follows that any solution of~Eq.~(\ref{cj}) can be represented as
a sum
$$
\sum_\lambda c_\lambda s_\lambda e^{f(\lambda)\beta},
$$
over all partitions~$\lambda$, for some coefficients~$c_\lambda$.
For the solution~$e^H$, these coefficients can be computed from the
initial value $H|_{\beta=0}=p_1$, which yields 
\begin{equation}\label{eigen}
e^H=\sum_\lambda s_\lambda(1,0,0,\dots) s_\lambda e^{f(\lambda)\beta}.
\end{equation}
In~\cite{O}, this form of the function~$e^H$ was deduced
from the representation theory of symmetric groups.
Note that $c_\lambda=(\dim~R_\lambda)/|\lambda|!$, where
$|\lambda|=\lambda_1+\lambda_2+\dots$ and $R_\lambda$
is the irreducible representation of the symmetric group~$S_{|\lambda|}$
corresponding to the partition~$\lambda$.

Taking the logarithm we obtain a few first terms in the expansion of~$H$:
\begin{eqnarray*}
H&=&p_1+\frac14(e^\beta-2+e^{-\beta})p_1^2+\frac14(e^\beta-e^{-\beta})p_2
+\frac1{36}(e^{3\beta}-9e^\beta+16-9e^{-\beta}+e^{-3\beta})p_1^3\\
&&\phantom{p_1}+\frac1{12}(e^{3\beta}-3e^\beta+3e^{-\beta}-e^{-3\beta})p_1p_2
+\frac1{18}(e^{3\beta}-2+e^{-3\beta})p_3+\dots
\end{eqnarray*}

\subsection{The exponent of the generating function for Hurwitz numbers
as an element of the semi-infinite Grassmannian}

In Sec.~2.2, the fact that the function~$e^H$ is a~$\tau$-function
for the KP hierarchy was established by 
a reference to Okounkov's paper~\cite{O}.
Here we present a more direct argument.

For~$\beta=0$, the function~$e^H=e^{p_1}$ is the $\tau$-function corresponding
to the semi-infinite space with the basis $\varphi_j(z)=z^je^{z^{-1}}$.
In other words, it is the result of the application of the
operator $e^{z^{-1}}$ {\rm(}the exponent being understood as the operator
of multiplication by~$z^{-1}${\rm)}
to the wedge product $z\wedge z^2\wedge z^3\wedge\dots$
{\rm(}the vacuum vector{\rm)}.

Indeed, for this basis $\varphi_j$, the coefficient of~$s_\lambda$
is nothing but $s_\lambda(1,0,0,\dots)$ because of the defining
equation~(\ref{Schur}) and the relation $s_k(1,0,0,\dots)=1/k!$
for the one-part partitions.

The most natural approach
to description of the evolution of this semi-infinite space along~$\beta$ 
consists in multiplying the coefficient of~$z^i$ in each $\varphi_j$
by~$e^{i(i-1)\beta/2}$. In other words, we apply to this space the
exponent of the diagonal matrix with the diagonal elements $i(i-1)\beta/2$.
This procedure produces a well-defined basis in the $\beta$-evoluting space.
For this basis, however, Eq.~(\ref{tau}) cannot be applied directly.
To avoid this difficulty, one can normalize the basic vector~$\varphi_j$
by dividing it by $e^{j(j-1)\beta/2}$ thus making it start with $z^j$:
\begin{equation}\label{ph}
\varphi_j(z)=\sum_{i\le j}\frac1{(j-i)!}e^{\frac{i(i-1)-j(j-1)}2\beta}z^i.
\end{equation}
Now it is obvious that Eq.~(\ref{tau}) gives the coefficient 
of~$s_\lambda(1,0,0,\dots)s_\lambda$
equal to $e^{f(\lambda)\beta}$, as desired, and we obtain

\begin{theorem}
The function~$e^H$ is the $tau$-function for the KP hierarchy corresponding
to the semi-infinite subspace with the basis~{\rm(\ref{ph})}.
\end{theorem}

This yields an independent proof of the following 

\begin{corollary}
The function~$H$
satisfies the KP hierarchy, in particular, it is a solution
to the KP equation~{\rm(\ref{KP})}.
\end{corollary}

The KP hierarchy degenerates into the KdV hierarchy 
if we are looking for solutions independent of variables
with even indices (that is, the derivatives with respect
to these variables vanish identically). In the normalization
of the present paper, the variables~$t_i$ are related to
the variables~$p_j$ by $t_i=\frac1{(2i-1)!!}p_{2i+1}$
which, in particular, makes~$F$ independent of $p_2,p_4,\dots$.
(Recall that $(2i-1)!!$ denotes the product of odd numbers
from 1 to $2i-1$, $(-1)!!=1$.)
Since the function~$F$ given by Eq.~(\ref{Wit})
satisfies the KdV hierarchy, its exponent $e^F$
is a $\tau$-function for the KdV
hierarchy. Its explicit presentation similar to that of~$e^H$
is given in~\cite{KS}. It would be interesting to deduce
the whole KdV-hierarchy for~$F$ from the KP-hierarchy for~$H$ directly,
without referring to the string and the dilaton equations.

\end{document}